\documentclass[11pt]{article}
\usepackage[dvips]{graphicx}
\usepackage{a4wide}
\usepackage{wasysym}
\usepackage[english]{babel}
\usepackage{amsfonts}
\usepackage{amssymb}
\newtheorem{defi}{\textbf{Definition}}[section]
\newtheorem{theo}[defi]{\textbf{Theorem}}
\newtheorem{lemma}[defi]{\textbf{Lemma}}
\newtheorem{prop}[defi]{\textbf{Proposition}}
\newtheorem{coro}[defi]{\textbf{Corollary}}

\newcommand{\dist}{{\rm dist}}
\newcommand{\Ang}{{\rm Ang}}
\newcommand{\MaxAng}{{\rm MaxAng}}
\def\proof{\noindent \textit{ Proof. }}

\title{Symbolic dynamics and relatively hyperbolic groups}

\author{Fran\c{c}ois Dahmani\footnote{Universit\'e Louis Pasteur,
Strasbourg, e-mail : dahmani@math.u-strasbg.fr}, Asl\i
Yaman\footnote{University of Southampton, e-mail :
a.yaman@maths.soton.ac.uk}}

\date{}

\begin{document}

\maketitle

\begin{center}

\footnotesize \textbf{Abstract :} We study the action of a
relatively hyperbolic group on its boundary, by methods of
symbolic dynamics. Under a condition on the parabolic
 subgroups, we show that this dynamical system is finitely
presented. We give examples where this condition is
satisfied, including geometrically finite
kleinian groups.

\normalsize

\end{center}

Associated to any word-hyperbolic group $\Gamma$, there is a
dynamical system arising from the action of $\Gamma$ on its Gromov
boundary $\partial \Gamma$. Already in \cite{Gro}, M.Gromov uses
methods of symbolic dynamics for the study of this action, and in
\cite{CP} (see also \cite{CP2}), M.Coornaert,
and A.Papadopoulos explain a way to
factorize such a dynamical system through a subshift of finite
type. They describe a finite alphabet $\mathcal{A}$, and a
subshift $\Phi \subset \mathcal{A}^{\Gamma}$, and they construct a
continuous equivariant, surjective map $\Phi \to \partial \Gamma$,
which encodes the action of $\Gamma$ on its boundary by a subshift
of finite type.

The action of a group $\Gamma$ on a compact metric space, $K$ is
expansive if there exists $\varepsilon>0$ such that any pair of
distinct points in $K$ can be taken at distance at least
$\varepsilon$ from each other by an element of $\Gamma$. It is
well known that the action of a hyperbolic group $\Gamma$ on
$\partial \Gamma$ is expansive. This property, together with the
existence of the coding given in \cite{CP}, makes the action of a
hyperbolic group, $\Gamma$, on its boundary, $\partial \Gamma$,
\emph{finitely presented} (see \cite{Gro}, \cite{CP}). In
\cite{Gro}, M.Gromov describes consequences of such a
presentation, like the rationality of some counting functions.

The aim of this paper is to state and prove similar properties for relatively 
hyperbolic groups, where parabolic subgroups are allowed. 
In general, in presence of parabolics, the study of dynamical 
properties becomes significantly more complicated.

After an idea of Gromov in \cite{Gro}, B.Farb \cite{F} and
B.Bowditch \cite{Brel} developed the theory of relatively
hyperbolic groups, as a generalization of geometrically finite
Kleinian groups. We will use for this work the definition of
relatively hyperbolic groups given by Bowditch in \cite{Brel}. A
group $\Gamma$ is hyperbolic relative to a family, $\mathcal{G}$,
of finitely generated subgroups of $\Gamma$ if it acts on an
hyperbolic fine graph, with finite stabilizers of edges, finitely
many orbits of edges, and such that the stabilizers of infinite
valence vertices are exactly the elements of $\mathcal{G}$ (see
Definition \ref{def;relhyp}). In \cite{F}, this definition is
equivalent to ``relatively hyperbolic with the property BCP'' (see the
appendix of \cite{theseDah}).

If one replaces ``fine'' by ``locally finite'' in above
definition, then $\mathcal{G}$ is empty and the group is
hyperbolic. In \cite{F}, Farb proves that the fundamental group of
a finite volume manifold of pinched negative curvature, with
finitely many cusps is hyperbolic relative to the conjugates of
the fundamental groups of the cusps, which are virtually
nilpotent. Sela's limit groups, or, finitely generated
$\omega$-residually-free groups are hyperbolic relative to their
maximal abelian non-cyclic subgroups, as shown in \cite{Dcomb}.

Bowditch describes a boundary for a relatively hyperbolic group in
\cite{Brel}. The group acts on this compact set as a convergence group, and the
elements of the family $\mathcal{G}$ are parabolic subgroups for this action.
Despite of those parabolic subgroups, the action is expansive
(Proposition \ref{prop;exp}). Although the construction of the
subshift of finite type given by Coornaert and Papadopoulos will
not work properly here (either one would need an infinite
alphabet, or the map $\Phi \to \partial \Gamma$ would not be well
defined) we found an intrinsic property of the maximal parabolic
subgroups that allows successful modifications.

An infinite group has its one-point compactification $G
\cup \{\infty\}$ finitely presented with special symbol if
there exists an alphabet $\mathcal{A}$, a subshift of finite type $\Phi \subset
\mathcal{A}^G$,  a continuous surjective $G$-invariant map $\Pi:\Phi \to (G \cup
\{\infty\})$, and a special symbol $\$ \in \mathcal{A}$, such that, for  $\sigma
\in \Phi$, $\Pi(\sigma) = g \in G$ if and only if $\sigma(g) = \$ $.
One can reformulate this definition including finite groups by considering their Alexandrof compactifications. If the action of a dicrete group on its Alexandrof compactification is finitely presented with special symbol then we say that the group has \emph{the property of special symbol}.

\begin{theo}
Let $(\Gamma, \mathcal{G})$ be a relatively hyperbolic group, and
$\partial \Gamma$ be its boundary (in the sense of Bowditch \cite{Brel}).

If each $G \in \mathcal{G}$ has the property of special symbol, then the action of $\Gamma$ on its
boundary $\partial \Gamma$ is finitely presented.

\end{theo}

One should remark that in \cite{CP} the subshift coding the action of hyperbolic groups to its boundary enables algorithmic computations. Here in order to simplify the proof we give a weak version of the result for relatively hyperbolic groups. However we refer to the thesis of second author (\cite{theseAs}) for a complete proof of the above theorem which enables computations.

We see that there is the natural question arising from above result. 

\medskip

\noindent {\bf Problem.} {\it Which groups have the property of special symbol?}

\medskip

 The following results answer partially this question and give some examples. However we do not have a complete list of groups satisfying the special symbol property.

\begin{theo}

If a group has the property of special symbol, then it is finitely generated.

\end{theo}

 A group $\Gamma$ is said
to be \emph{poly-hyperbolic} if  there is a sequence of subgroups $\{1\}=N_0
\lhd N_1 \lhd ... \lhd N_{k-1} \lhd N_k = \Gamma$, with all the quotients
$N_{i+1}/N_i$ hyperbolic.

\begin{theo}
 Poly-hyperbolic groups have the property of special symbol. 
 In particular, this includes poly-cyclic groups.
\end{theo}

\begin{coro}

A geometrically finite group \cite{Bgf} of isometries of a simply connected 
manifold of pinched negative curvature acts as a finitely presented system on
its limit set. 

 \end{coro}

\hskip .4cm

We give in section 1 definitions related to symbolic dynamics. In
section 2 we define relatively hyperbolic groups, their boundaries and
introduce some tools such as "angles" and "cones". We prove the
Theorem 0.1 in Section 3. The subshift we construct will produce
objects which are local Busemann functions on the fine hyperbolic
graph associated to the relatively hyperbolic group. To associate
a point in the boundary to an element of the subshift, we consider
its gradient lines. We prove that they converge to points at
infinity, and we make sure that, for a given element of the
subshift, all the gradient lines converge to the same point. For
this we use the property of special symbol for each stabilizer of
infinite valence vertex. In Section 4 we study this property of
special symbol, and in particular, prove Theorem 0.2.

We want to thank B.Bowditch, and T.Delzant, for their many comments, and
M.Coornaert for his useful explanations. We also thank the referee, for
his/her useful comments.

    \section{Definitions, symbolic dynamics}

We borrow the following definitions (1.1 to 1.4) from Gromov \cite{Gro} 8.4, and
Coornaert and Papadopoulos \cite{CP}, chapter 2. See also Fried \cite{Fri}.

\begin{defi}(Shift, subshift, subshift of finite type)
\label{def;through}(\cite{Gro}, 8.4, \cite{CP}, Chap. 2)

If $\mathcal{A}$ is a finite discrete alphabet and $\Gamma$ is a group,
$\mathcal{A}^{\Gamma}$, with the product topology, is the total
shift of $\Gamma$ on $\mathcal{A}$. It admits a natural left
$\Gamma$-action given  by  $(\gamma\sigma)(g)= \sigma(\gamma^{-1}g)$
for all $g\in \Gamma$, $\sigma \in \mathcal{A}^{\Gamma}$.

A closed $\Gamma$-invariant subset of $\mathcal{A}^{\Gamma}$ is
called a subshift.

A cylinder $\mathcal{C}$ is a subset of the total shift such that
there exists a finite set $F\subset \Gamma$, and a family of maps
$M\subset \mathcal{A}^{F}$ with $$\mathcal{C} = \{\sigma \in
\mathcal{A}^{\Gamma} \, | \; \sigma|_{F} \in M \}.$$

$\Phi$ is a \emph{subshift of finite type} if there exists a
cylinder $\mathcal{C}$ such that $\Phi = \bigcap_{\gamma \in
\Gamma} \gamma^{-1}\mathcal{C}$.

\end{defi}

The subshifts of finite type are subshifts, but the
cylinders are not $\Gamma$-invariant.

\begin{defi}(Dynamical systems of finite type) \label{def;def2}\cite{Gro},  \cite{CP}

Let $\Gamma$ act on a compact set $K$. The dynamical system is of
\emph{finite type} if there exists a finite alphabet $\mathcal{A}$, a
subshift of finite type $\Phi \subset \mathcal{A}^{\Gamma}$ and a
continuous, surjective, $\Gamma$-equivariant map $\pi : \Phi \to
K$.

\end{defi}

\textit{Example :} Set $\Gamma =\mathbb{Z}$, and
$\mathcal{A}=\{a,b,\$\}$. The cylinder $\mathcal{C}$ is the set of
the maps that agree on $F=\{0,1\}$ with one of the maps in $M=\{
m_i,\; i=1...4\}$ where $(m_1: 0,1 \to a,a)$, $(m_2:0,1 \to a,\$) $,
$(m_3: 0,1 \to \$,b) $ and $(m_4 : 0,1 \to b,b) $.

Let $\Phi$ be the subshift of finite type defined by the cylinder
$\mathcal{C}$, i.e $\Phi = \bigcap_{n \in \mathbb{Z}}
n+\mathcal{C}$, where $n+\mathcal{C}$ is the translate of the cylinder
$\mathcal{C}$ of the total shift $\mathcal{A}^{\mathbb{Z}}$,
by the element $n$ of $\mathbb{Z}$ .
One can check that the elements of $\Phi$ are the constant word on $a$, the
constant word on $b$ and all the words $(...aaa\$bbb...)$ beginning by $a$,
until there is a $\$ $ on the $n^{th}$ letter ($n\in \mathbb{Z}$) and then $b$.
Although for this example $\Phi$ is countable, in general subshifts of finite
type are not countable.

Now consider the compact set $K=\mathbb{Z}\cup\{\infty\}$ with the
usual topology. There is a natural left action of $\mathbb{Z}$ on
$K$, fixing the infinity point. Consider the map $\pi : \Phi \to
K$ that sends $(...aaaa...)$ and $(...bbbb...)$ on $\infty$, and
$(...aaa\$bbb...)$ on $n\in \mathbb{Z}$ where $n$ is the index of
the letter $\$ $. The map $\pi$ is surjective, continuous and
equivariant, and therefore the action of $\mathbb{Z}$ on $K$ is of
finite type.

We now continue with definitions. One can refine the property of
being a dynamical system of finite type with the following.

\begin{defi}(Expansivity)

The action of a group $\Gamma$ on a compact set $K$ is
\emph{expansive} if there exists $U$ a neighbourhood  of the
diagonal $\Delta\subset K\times K$  such that  $\Delta =
\bigcap_{\gamma \in \Gamma} \gamma U$.

\end{defi}

Note that, if the compact set is metric, this is equivalent to the
definition of expansivity given in introduction (see
\cite{CP}, Proposition 2.3).

\begin{defi}(Finitely presented dynamical systems)\cite{Gro}, \cite{CP}

Let $\Gamma$ act on a compact $K$. The dynamical system is
finitely presented if it is both of finite type and expansive.

\end{defi}

If one has a subshift of finite type $\Phi \subset
\mathcal{A}^{\Gamma}$ and a surjective continuous equivariant map
$\pi : \Phi \to K$, the expansivity of the action of $\Gamma$ on
$K$ turns out to be equivalent to the fact that the subshift $\Psi
\subset (\mathcal{A}\times \mathcal{A})^{\Gamma}$ defined by
$[(\sigma_1 \times \sigma_2) \in \Psi] \Leftrightarrow
[\pi(\sigma_1)=\pi(\sigma_2)]$, is of finite type (cf \cite{CP} chapter 2).

If $\Gamma$ is a infinite discrete group, it acts on its one-point
compactification $\Gamma \cup \{\infty\}$ by multiplication of the  left, hence
fixing the point at infinity. This one-point compactification is also the
Alexandrov compactification of $\Gamma$, but if $\Gamma$ is finite, its Alexandrov
compactification is itself.

\begin{defi}(The property of special symbol) \label{def;fp+sc}

A discrete group $\Gamma$ is has
\emph{the property of special symbol} if the $\Gamma$-action on the Alexandrov compactification, $\Gamma \cup \{\infty\}$ (or respectively on $\Gamma$ if $\Gamma$ is finite), is finitely presented by a subshift $\Phi\subset \mathcal{A}^{\Gamma}$ and if the presentation map $\pi : \Phi \to \Gamma \cup \{\infty\}$ (resp. $\pi : \Phi \to \Gamma$) satisfies

$$ \exists \$ \in \mathcal{A} \; \quad (\pi(\sigma)=\gamma \in
\Gamma) \Leftrightarrow (\sigma(\gamma)=\$) $$
\end{defi}

Note that in this case, the property of expansivity of $\Gamma$, on
$K=\Gamma\cup\{\infty\}$ is always satisfied (consider $U =
(\Gamma \times \Gamma) \cup  (\{\infty\}\times
(K \setminus \{1_\Gamma\})) \cup ((K \setminus \{1_\Gamma\}) \times
\{\infty\})$ in the definition 1.3 ). The example of dynamical system of finite type described
previously has the property of special symbol of $\mathbb{Z}$. Note
also that all finite groups have trivially the property of special symbol.

We give in section \ref{examples} several examples of groups which has the property of special symbol.

   \section{About Relatively Hyperbolic Groups}

\subsection{Definitions}

A \emph{graph} is a set of \emph{vertices} with a set of
\emph{edges}, where an edge is an unordered pair of vertices. One can equip the
geometrical realization of a graph with a metric where edges have length 1. Thus 
this geometrical realization allows to consider simplicial, geodesic,
quasi-geodesic and locally geodesic paths in a graph. (see for example \cite {Bgrhyp})

\begin{defi}(Circuits)

A \emph{circuit} in a graph is a simple simplicial loop, i.e
without self intersection.

\end{defi}

In \cite{Brel}, B. Bowditch introduces the notion of fineness of a
graph.

\begin{defi}(Fineness)\cite{Brel}

A graph $\mathcal{K}$ is \emph{fine} if for every $L>0$, and for every
edge $e$, the set of the circuits of length less than $L$,
containing $e$ is finite. It is \emph{uniformly fine} if  this set
has cardinality bounded above by a constant depending only on $L$.

\end{defi}

\begin{defi}(Relatively Hyperbolic Groups)\label{def;relhyp}\cite{Brel}

A group $\Gamma$ is \emph{hyperbolic relative} to a family of
subgroups $\mathcal{G}$, if it acts on a Gromov-hyperbolic, fine
graph $\mathcal{K}$, such that stabilizers of edges are finite,
 such that there are finitely many orbits of edges, and such that the
stabilizers of the vertices of infinite valence are finitely generated and they are exactly the elements of $\mathcal{G}$.

\end{defi}

With an abuse of language, we will say that the pair $(\Gamma,
\mathcal{G})$ is a relatively hyperbolic group, and that
$\mathcal{K} $ is a  graph \emph{associated} to it.

We note that as there are finitely many orbits of edges, a graph
associated to a relatively hyperbolic group is uniformly fine.
Note also that the graph $\mathcal{K}$ associated to $(\Gamma,\mathcal{G})$,
can be chosen to be without global cut point, and  with positive hyperbolicity
constant $\delta$. This will ensure that all the angles (introduced in the
following arguments) are finite.

 \subsection{Angles}

For any graph, one can define a notion of \emph{angle} as follow.

\begin{defi}(Angles)

Let $\mathcal{K}$ be a graph, and let $e_1=(v,v_1)$ and
$e_2=(v,v_2)$ be edges with one common vertex $v$. The angle
$\Ang_v(e_1,e_2)$, is the length of the shortest path from $v_1$
to $v_2$, in $\mathcal{K}\setminus \{v\}$ ($+\infty$ if there is
none).

\end{defi}

The angle $\Ang_v(p,p')$ between two simple simplicial (oriented)
paths $p$ and $p'$, starting from a common vertex $v$, is the angle between
their first edges after this vertex.

If $p$ is a simple simplicial path, and $v$ one of its vertices,
$\Ang_v(p)$ is the angle between the consecutive edges of $p$ at
$v$, and its maximal angle $\MaxAng(p)$ is the maximal angle
between consecutive edges of $p$.

In the notation $\Ang_v(p,p')$, we will sometimes omit the
subscript if there is no ambiguity.

\begin{prop}(Some useful remarks)\label{prop;remarks}

 1.   $\Ang_v(e_1,e_3)\leq \Ang_v(e_1,e_2)+\Ang_v(e_2,e_3)$
 when $e_i$ are edges incident on a vertex $v$.

 2.   If $\gamma$ is an isometry of $\mathcal{K}$, $\Ang_v(e_1,e_2)=
\Ang_{\gamma v}  (\gamma.e_1,\gamma.e_2)$.

 3.  Any circuit of length $L \geq 2$ has a maximal angle
 less than $L-2$.

\end{prop}

\proof The first remark is the triangular inequality for the length distance
of $\mathcal{K}\setminus \{v\}$. The second statement is obvious.
Finally, if $e_1=(v_1,v)$ and $e_2=(v,v_2)$ are two consecutive
edges in the circuit, the circuit itself gives a path of length
$L-2$ from  $v_1$ and $v_2$ avoiding $v$. $\square\smallskip$

Here is an important property of angles.

\begin{lemma}(Large angles in triangles)\label{lem;angles_triangles}

Let $[x,y]$ and $[x,z]$ be geodesic segments in a
$\delta$-hyperbolic graph, and assume that
$\Ang_x([x,y],[x,z])=\theta \geq 50\delta$. Then the concatenation
of the two segments is still a geodesic. Moreover any geodesic
segment $[y,z]$ will contain $x$ and $\Ang_x([y,z]) \geq \theta
-50\delta$.

\end{lemma}

\proof Let $[y,z]$ be a geodesic, defining a triangle $(x,y,z)$, which is
$\delta$-thin, that is : any segment $[y,z]$ is in the $\delta$-neighbourhood of
the set $[x,y] \cup [x,z]$.  We consider the vertices $y'$ and $z'$ on $[x,y]$
and $[x,z]$ located at distance $10 \delta$ from $x$. They are not
$3\delta$-close to each other. Indeed, if they were, there would be a loop of
length less than $23\delta$ containing $x$ and the firsts edges of $[x,y]$ and
$[x,z]$, and not returning to $x$, which contradicts the fact that the angle of
these path at $x$ is more than $50\delta$. Therefore, they are
$\delta$-close to the segment $[y,z]$, and we set $y''$ and $z''$ the
corresponding points on $[y,z]$. This gives a loop of length less than
$(2\times 10\delta+2\delta)\times 2  \le 50\delta$, containing $x$, consisting
of $[x,y']$, $[y',y'']$, $[y'',z'']$, $[z'',z']$, and $[z',x]$. As the small
segments $[y',y'']$ and $[z'',z']$ are $10\delta$ far away from $x$, they do
not contain $x$. Thus the third property of Proposition \ref{prop;remarks}
proves that $x\in [y'',z'']$, and $\Ang_x([y'',z'']) \geq \theta -50\delta $,
and therefore $\Ang_x([y,z]) \geq \theta-50\delta$ . $\square\smallskip$

\subsection{Cones}

\begin{defi}(Cones)

Let $\mathcal{K}$ be a graph, let $d$ and $\theta$ be positive
numbers.  The \emph{cone} centred at an edge $e=(v,v')$, of
radius $d$ and angle $\theta$ is the set of vertices $w$ at distance
at most $d$ from $v$ and such that there exists a geodesic segment
$[v,w]$ the maximal angle and the angle with $e$ of which are
less than $\theta$, i.e. :

$$ Cone_{d,\theta}(e,v)  =   \{w\in \mathcal{K} \,| \,  \dist(w,v)\!\leq\!
d,\,\MaxAng[v, w]\!\leq\!\theta, \Ang_v(e,[v, w])\! \leq\! \theta
\} $$

\end{defi}

\begin{prop}(Bounded angles imply local finiteness)

Let $\mathcal{K}$ be a fine graph. Given an edge $e$ and
$\theta>0$, there exists only finitely many edges $e'$ such that
$e$ and $e'$ have a common vertex, and $\Ang(e,e')\leq \theta$.

\end{prop}

\proof Only finitely many circuits shorter than $\theta$ contain $e$.
$\square\smallskip$

\begin{coro}(Cones are finite)\label{cor;conesarefinite}

In a fine graph, the cones are finite sets of vertices. If the
graph is uniformly fine, the cardinality of $Cone_{d,\theta}(e,v)$
can be bounded above by a function of $d$ and $\theta$.
\end{coro}

\proof Consider a cone $Cone_{d,\theta}(e,v)$. We argue by induction on $d$. If $d=1$,
the result is given by the previous proposition. If $d>1$, we remark that
$Cone_{d,\theta}(e,v)$ is contained in the union of cones of angle $\theta$ and
radius $1$, centered at edges whose vertices are both in
$Cone_{(d-1),\theta}(e,v)$. If the latter is finite, the union is also finite.
$\square\smallskip$

Here is a visibility property for hyperbolic fine graphs. It is an usual result
for proper hyperbolic graphs.

\begin{prop}\label{prop;visual}(Visibility property in fine hyperbolic graphs)

 Let $\mathcal{K}$ a hyperbolic fine graph, and $\partial \mathcal{K}$ its
Gromov boundary. Then for all $\xi \in \partial \mathcal{K}$, and for all vertex
$v$ in $\mathcal{K}$, there exists a geodesic ray $\rho = [v,
\xi)$ in $\mathcal{K}$.

For all $\xi$ and $\xi'$  in $\partial \mathcal{K}$, there exists a bi-geodesic
ray $\rho' = (\xi, \xi')$ in $\mathcal{K}$.
\end{prop}

\proof Let us prove the first assertion. Let $(v_n)$ be a sequence of vertices
converging to $\xi$ in the sense of the Gromov topology. Consider two segments
$[v, v_n]$, and $[v,v_m]$. If their angle at the vertex $v$ is greater than
$50\delta$, we can apply Lemma \ref{lem;angles_triangles}  to deduce that the
gromov product $(v_n \cdot v_m)_v$ is equal to zero. Therefore, there exists
$n_0$ such that for all $n\geq n_0$, $Ang_v([v, v_n],[v, v_{n_0}] \le 50\delta$.
Therefore, the first edges of these segments are all in a cone, and as the graph
is fine, it is a finite set. After extraction of a subsequence, we can assume
that these edges are all equal. The diagonal extraction process gives a
subsequence $(v_{\sigma(n)})_n$ such that $[v, v_{\sigma(n)}]$ coincide with
$[v, v_{\sigma(m)}]$ for all $m\geq n$, on a subsegment $s_n$ of length $n$. As
$s_n$ is a subsegment of $s_m$ for all $m\ge n$, their union is a geodesic ray
$[v, \xi)$.

For the second assertion, we choose two rays $[v,\xi)$ and $[v,\xi')$. Let
$(v_n)_n$ be a sequence of vertices that converges to $\xi$, on  $[v,\xi)$
and $(v'_n)_n$ another sequence that converges to $\xi'$ on  $[v,\xi')$.
Let $d =(\xi \cdot \xi')_v +100\delta$, and let $r$ and $r'$ the subsegments of
the rays $[v,\xi)$ and $[v,\xi')$, of length $d$.  For all $n$ and $m$
sufficiently large, the geodesic segments $[v_n,v_m]$ intersect the cone
of radius $d$ and of angle $\MaxAng (r) + \MaxAng (r') + Ang_v(r,r') +100\delta$
centered at the first edge of $r$.

Therefore, as the cone is finite, one can find a subsequence of $v_n$ and
$v'_n$ and a point $p$ in the cone, such that every segment $[v_n, v'_m]$
contains $p$. Now the diagonal process of the proof of the first assertion gives
a sequence of segment converging to a bi-geodesic. $\square\smallskip$

We will use the following theorem, which is a reformulation of a
result in \cite{Dah}.

\begin{theo}\label{thm;Rips}

Let $\Gamma$ be a relatively hyperbolic group and $\mathcal{K}$ be
an associated graph, which is $\delta$-hyperbolic. There exists an
aspherical (in particular simply connected) simplicial complex such
that its vertex set is the one
of $\mathcal{K}$, and such that each simplex has all its vertices
in a same cone of $\mathcal{K}$, of radius $10\delta +10$, and
angle $100\delta +30$.
\end{theo}

 \proof In \cite{Dah}, the first author defines the \emph{relative Rips
complex} $P_{d,r}(\mathcal{K})$ for a relatively hyperbolic group.
It is the maximal complex on the set of vertices of $\mathcal{K}$
such that an edge is between two vertices if, in $\mathcal{K}$, a
geodesic of length less than $d$ and maximal angle less than $r$
links them. Although in \cite{Dah}, the notion of angle is
replaced by ``length of traveling in cosets'', the proof of
Theorem 6.2 remains the same, and gives the asphericity of
$P_{d,r}(\mathcal{K})$ for large $d$ and $r$. Theorem
\ref{thm;Rips} follows. $\square\smallskip$

   \subsection{Boundary of a relatively hyperbolic group}

Let $(\Gamma,\mathcal{G})$ be a relatively hyperbolic group, and
let $\mathcal{K}$ be an associated graph. In \cite{Brel}, Bowditch
defines the (dynamical) boundary $\partial \Gamma$ of $\Gamma$ by
$\partial \Gamma = \partial \mathcal{K} \cup \mathcal{V}_{\infty}$
where $\partial \mathcal{K}$ is the Gromov boundary of the
hyperbolic graph $\mathcal{K}$, and $\mathcal{V}_{\infty}$ is the
set of vertices of infinite valence in $\mathcal{K}$. This
boundary admits a natural topology of metrisable compact set (see
\cite{Brel}, \cite{Bbord}, \cite{Asli}). Let us recall a convergence
criterion for this topology. Let $(v_n)_n$ be a sequence of
vertices of infinite valence. If there is a point $\xi$ in $\partial
\mathcal{K}$ such that the sequence of Gromov products $(v_n \cdot \xi)_{v_0}$
tends to infinity, then $(v_n)_n$ converges to $\xi$ for the topology of
$\partial \Gamma$. If there is a vertex $v'$ in $\mathcal{K}$ such that every
geodesic segments $[v_0, v_n]$ contains $v'$ and contains an edge $e_n = (v',
v'_n)$ such that all the vertices $v'_n, \, n \in \mathbb{N}$ are distinct, then
the sequence $(v_n)_n$ converges to $v'$ in the topology of $\partial \Gamma$.
For a sequence $(\xi_n)_n$ of points in $\partial \mathcal{K}$, the conditions
are similar : one needs only to change the segments $[v_0,v_n]$ into geodesic
rays $[v_0,\xi_n]$.

\section{Finite presentation of the boundaries of a relatively hyperbolic group.}

We will prove the following theorem.

\begin{theo}\label{theo;main}

Let $(\Gamma,\mathcal{G})$ be a relatively hyperbolic group. If,
for each $G \in \mathcal{G}$, the action of $G$ on its Alexandrov
one-point compactification $G \cup \{\infty\}$ has the property of special symbol, then the action of $\Gamma$ on its
 boundary $\partial \Gamma$ is finitely presented.
\end{theo}

Let $(G_i)_{i=1..m}$ be a finite family of representatives of conjugacy
 classes of parabolic subgroups
in $\Gamma$, each stabilizing an infinite valence vertex $p_i$ in some hyperbolic fine graph $\mathcal{K}$. 
We also choose, for each $i\leq m$, an arbitrary edge $e_i$ in $\mathcal{K}$, adjacent to
$p_i$. 
We assume that each $G_i$ acts on $G_i\cup
\{\infty\}$ as a finitely presented system with a special symbol. That means
that we have an alphabet $\mathcal{A}_i$, a finite subset
$F_i\subset G_i$ and a set $M_i$ of maps from $F_i$ to
$\mathcal{A}_i$ which define a cylinder, hence a subshift of
finite type $\Phi_i$. We will denote by $\$ $ the special symbol in $\mathcal{A}_i$,
without distinguishing the indices $i$.

Without loss of generality, we can choose that the graph $\mathcal{K}$ is such that for all $i$, 
and for all $\gamma \in F_i$, $\Ang(e_i,\gamma e_i) \leq 1$. We choose also $\mathcal{K}$ to be without cut point (all the angles are finite).

\subsection{Busemann and radial cocycles}

\begin{defi}(Busemann function)(see \cite{Gro} 7.5.C, and \cite{CP} chap. 3, section 3)

Let $\rho : [0,\infty) \to \mathcal{K}$ be a geodesic ray starting
at $v_0$. The Busemann function $h_{\rho} :
\mathcal{V} \to \mathbb{Z}$ of $\rho$ is defined by the limit (which always
exists and is finite) $h_{\rho}(v) = \lim_{n \to \infty} (\dist(v,\rho(n)) -
n)$.

\end{defi}

\begin{defi}(Busemann cocycles)(\cite{Gro} 7.5.E, \cite{CP} chap. 3)

Let $h_{\rho}$ be a Busemann function. The cocycle associated to
$h_{\rho}$ is $\varphi_{\rho} : \mathcal{V}\times \mathcal{V} \to
\mathbb{Z}$ defined by $\varphi_{\rho}(w,v) = h_{\rho}(v) -
h_{\rho}(w)$. A gradient line of $\varphi_{\rho}$ is a sequence of
vertices $(v_n)_n$ such that $\varphi_{\rho}(v_{i+1},v_i)=1$ for
all $i$.

\end{defi}

The proof of the next lemma can be found in \cite{CP} (Proposition 4.2), for locally finite graphs, which is not the case here. 
But in fact, the proof given depends only a visibility property of the boundary, 
which is satisfied here (Proposition \ref{prop;visual}).

\begin{lemma}\label{lem;gradlines}

If $\varphi$ is a Busemann cocycle associated to $\rho$, then a
gradient line is  a sequence of vertices of a geodesic ray
asymptotic to $\rho$. Moreover there is a gradient line starting
from each vertex.

\end{lemma}

\begin{defi}(Radial cocycles)

Let $p$ be a vertex of infinite valence in $\mathcal{K}$. The radial
cocycle associated to $p$ is $\varphi_{p} : \mathcal{V}\times
\mathcal{V} \to \mathbb{Z}$ defined by $\varphi_h(w,v) =
\dist(v,p) - \dist(w,p)$. A gradient line of $\varphi_{p}$ is a
finite family of vertices $(v_n)_{0\leq n \leq m}$ such that
$\varphi_{p}(v_{i+1},v_i)=1$ for all $i$, and $v_m = p$.

\end{defi}

The next lemma is direct by definition.

\begin{lemma}\label{lem;obvious}

If $\varphi$ is a radial cocycle associated to ${p}$, then its
different gradient lines are exactly the sequences of vertices of
geodesic segments ending at $p$.

\end{lemma}

Let $\delta$ be a positive hyperbolicity constant of $\mathcal{K}$.
We set $\theta \geq 2000 \delta$ and such that for every vertex of finite valence in $\mathcal{K}$, for every pair of edges adjacent to this vertex, their angle is at most $\theta$.

\begin{prop}(Properties of Busemann and radial cocycles)\label{prop;locbusrad}

Let $\varphi$ be a Busemann or a radial cocycle. Then :

\begin{enumerate}
  \item (Integral values) For $x,y$ adjacent vertices, $\varphi(x,y)$ is $0$, $1$ or $-1$.
  \item (Cocycle) For all $x,y,z$ $\varphi(x,y)+\varphi(y,z)+\varphi(z,x)=0$.
  \item (Geodesic extension) Let $\xi$ be a point of $\partial \Gamma$
(\emph{i.e.} a point of $\partial \mathcal{K}$ or a vertex af infinite valence of
$\mathcal{K}$), and let $l=[v,\xi)$ be a gradient line and $[x,v]$ a geodesic
segment (of length and maximal angle less than $\theta$), such that
$\Ang_v([x,v]\cup[v,\xi)) \geq \theta$, then $[x,v]\cup[v,\xi)$ is a gradient
line.

  \item (Exits) If $v$ is a vertex of finite valence, then there exists $w$
adjacent to $v$ with $\varphi(w,v)=1$.

\end{enumerate}

\end{prop}

\proof Properties 1 and 2 are obvious. Property 4 is consequence of the
Lemmas \ref{lem;gradlines} and \ref{lem;obvious}. Property 3
deserves a proof here. We precise that the assumption of maximal angle for the segment $[x,v]$ is unnecessary here, 
but will be useful elsewhere in this paper, where we prove that some other objects satisfy the same properties 
(see \ref{lem;thirdpoint}). By Lemma \ref{lem;gradlines} (if $\xi \in
\partial \mathcal{K}$) and Lemma \ref{lem;obvious} (if $\xi$ is a
vertex of infinite valence), any gradient line from $x$ is a geodesic ray
$[x,\xi)$ and produces a triangle $(x,v,\xi)$, which, by
assumption, has a large angle at $v$. Hence, by Lemma
\ref{lem;angles_triangles}, any ray $[x,\xi)$ contains $v$.
$\square\smallskip$

\subsection{Shift and subshift}

For all $i=1\dots m$ we set $F'_i = \{ \gamma \in G_i \, | \; \Ang_{p_i}(e_i, \gamma e_i)
\leq \theta/2 \}$. This set contains $F_i$.

We fix a vertex $v_0$ and an edge $e_0=(v_0,v)$. We choose $R$ and
$\Theta$ sufficiently large, such that for all $i=1\dots m$,
$Cone_{10\theta ,10 \theta}(e_i, p_i) \subset Cone_{R,\Theta}(e_0,
v_0) $. 
 Note that $\theta \geq 2000 \delta$, hence $R$ and $\Theta$ are 
greater than $(100 \delta + 30)$, the constant given
 by Theorem \ref{thm;Rips}.

 Let $\mathcal{A}'$ denote the set of all possible restrictions of
Busemann and radial cocycles on $Cone_{R,\Theta}(e_0,v_0) \times
Cone_{R,\Theta}(e_0,v_0)$. We set $\mathcal{A}'' =\mathcal{A}_1
\times \dots \times \mathcal{A}_m$. We choose our alphabet to be
$\mathcal{A}=\mathcal{A}' \times \mathcal{A}''$

\begin{lemma}

$\mathcal{A}=\mathcal{A}'\times \mathcal{A}''$ is finite.

\end{lemma}

\proof Cones are finite, and cocycles have integral values bounded by the
diameter. $\square\smallskip$

An element $\psi$ of  $\mathcal{A}^\Gamma$ is a map from $\Gamma$
to $\mathcal{A}=\mathcal{A}' \times \mathcal{A}''$. Thus it has
coordinates $\psi_0 : \Gamma \to \mathcal{A}'$ and $\psi_i :
\Gamma \to \mathcal{A}_i$ for all $i$. Hence, $\psi_0(\gamma)$ is
a map from $Cone_{R,\Theta}(e_0,v_0) \times
Cone_{R,\Theta}(e_0,v_0)$ to $\mathbb{Z}$, whereas
$\psi_i(\gamma)$ is in $\mathcal{A}_i$ for $i\geq 1$.

Let $F$ be the set of elements in $\Gamma$ such that the vertices
of $\gamma.e_0$ are both in $Cone_{R,\Theta}(e_0,v_0)$. As
stabilizers of edges are finite, $F$ is a finite set.

Let $\mathcal{C}$ be the cylinder (in the sense of Definition
\ref{def;through}) defined on $F$ so that $\psi \in \mathcal{C}$
if the three next conditions, which concern only finitely many
elements of $\Gamma$, are fulfilled :

\begin{itemize}
  \item $[\psi_0(\gamma)](v_1,v_2)=[\psi_0(1_\Gamma)](\gamma^{-1}v_1,\gamma^{-1}v_2)$
whenever $v_1, v_2, \gamma^{-1}v_1, \gamma^{-1}v_2$ are all in
$Cone_{R,\Theta}(e_0,v_0)$.
  \item  $\psi_i |_{F_i}$ is in $M_i$.
  \item  for $\gamma \in F_i$, for $v$ such that $\gamma.e_i =(p_i,v)$,
and for $w$ such that $[w,p_i]$ is a geodesic segment of length,
and maximal angle, less than $\theta$,  containing $v$, one has
$[\psi_0(\gamma)](\gamma^{-1}w,p_i)\geq (1-\dist(w,p_i))$ only if
there exists $\gamma' \in F'_i$ such that $\psi_i(\gamma\gamma')=\$ $.

\end{itemize}

Denote by $\Phi$ the subshift of finite type  $\Phi =
\bigcap_{\gamma\in \Gamma} \gamma \mathcal{C}$.  

We will see in the next lemmas that the second point of the subshift allows us to say that at a vertex of infinite valence there is at most one edge issued from this vertex for which the $i$th coordinate of $\psi$ takes the symbol $\$ $. Moreover the third point of the subshift will ensure that the gradient lines escape a vertex of infinite valence only around the edge which takes the symbol $\$ $.

\begin{lemma}(Globalisation) \label{lem;global}

Let $\psi \in \Phi$. If $v$ and $v'$ are vertices in $\gamma
Cone_{R,\Theta}(e_0,v_0)$, we set $\varphi_{\psi} (v,v') = \psi_0
(\gamma) (\gamma^{-1}v,\gamma^{-1}v')$. Then the map
$\varphi_{\psi}$  is well defined. Moreover $\varphi_{\psi}$ can be extend to all pair of vertices of ${\mathcal K}$.

\end{lemma}

\proof Because of the first property of the definition of $\mathcal{C}$,
the formula given for $\varphi_{\psi}(v,v')$ does not depend on
the choice of possible $\gamma$, and therefore, the map is well
defined. Now, the map $\varphi_{\psi}$ is defined on pairs of vertices lying in
a same translate of $Cone_{R,\Theta}(e_0,v_0)$ and takes integral values, (first property of Proposition \ref{prop;locbusrad}) on the translate of this cone as it is a restriction of a busemann or a radial cocycle by definition of our alphabet. Thus it can be seen as a 1-cochain defined on the relative Rips complex given by Theorem \ref{thm;Rips}, which is simply connected. As it is a
cocycle, it is a coboundary, and there is a map $\tilde{\varphi}$
defined on the set of vertices of $\mathcal{K}$ such that
$\varphi_\psi (w,v)= \tilde{\varphi}(w)- \tilde{\varphi}(v)$ for
all $v,w$ lying in a translate of $Cone_{R,\Theta}(e_0,v_0)$. This
formula allows to extend the cocycle $\varphi_\psi$ to all pair of vertices
(not only those in a same cone of radius $R$ and angle $\Theta$).$\square\smallskip$

>From now on we will assume that  $\varphi_{\psi}$ is defined globally.

\begin{lemma}(About the $\psi_i$, $i\geq 1$)

Let $\psi \in \Phi$, and  $\gamma \in \Gamma$,  for all $i$,
$\psi_i |_{\gamma.G_i}$ is an element of $\Phi_i$.

\end{lemma}

\proof By definition of $\mathcal{C}$, for all $\gamma$, and all $g_i\in
G_i$, $\psi_i |_{\gamma g_i F_i}$ is in $M_i$. $\square\smallskip$

\begin{lemma}(About $\psi_0$)\label{lem;thirdpoint}

For all $\psi \in \Phi$, $\varphi_{\psi}$ satisfies each property of
Proposition \ref{prop;locbusrad} and that $|\varphi_{\psi} (w,v)| \leq \dist(w,v)$, for all $v$ and $w$.

\end{lemma}

\proof Given a line $l$ connecting $v$ and $w$, we have $\varphi_{\psi}(v_{n+1},v_n)\leq 1$ for every consecutive vertices $v_{n}$ and $v_{n+1}$ on $l$. Therefore, as
$\varphi_{\psi}$ is a global cocycle, $\varphi_{\psi}(v_{m},v_n)\leq m-n$, which proves the claim. Properties 1 and 4 of Proposition \ref{prop;locbusrad}
are satisfied because each element of our alphabet satisfy them in
a cone. Property 2 (cocycle property) fallows directly from the globalisation of $\varphi_{\psi}$ by \ref{lem;global}. For property 3, we notice that for $\gamma_0$ arbitrary, $\$ $ can
appear at most once in the set of values $\psi_i(\gamma_0 \gamma)$
as $\gamma$ ranges over $G_i$. If $[v,\xi)$ is a gradient line,
with $v = \gamma_0 p_i$, for some $i$ and some $\gamma_0$, then,
for $\gamma$ such that $\psi_i(\gamma_0 \gamma)=\$ $, there exists
$\gamma'$ in $F'_i$ such that $(\gamma_0\gamma \gamma').e_i$ is
the first edge of $[v,\xi)$. If $[v,x]$ is a segment such that
$\Ang_v([v,x],[v,\xi))\geq \theta$, then, for all $\gamma'' \in
F'_i$, the edge $(\gamma_0\gamma \gamma'').e_i$ is not on $[v,x]$.
Therefore, by the third property of the definition of the
cylinder, $\varphi_{\psi}(x,v)\leq - \dist(x,p_i)$. And this,
together with $|\varphi_{\psi}(x,v)|\leq \dist(x,v)$, gives
$\varphi_{\psi}(x,v)=- \dist(x,p_i)$ . In other words, $[x,v]\cup
[v,\xi)$ is a gradient line. $\square\smallskip$

\subsection{The presentation $\Pi : \Phi \to \partial \Gamma$}

Given an element of $\Phi$, we want to associate canonically an
element of $\partial \Gamma$.

\begin{defi}(Gradient lines)

Let $\psi$ be an element of $\Phi$. A gradient line $l_\psi$ of
$\psi$ is a finite or infinite sequence $(v_n)_{n\geq 0}$ of
vertices in $\mathcal{K}$ such that $\varphi_{\psi}(v_{n+1},v_n) =
1$ for all $n$. Moreover, it is finite only if for the last index
$m$, every neighbour $v$ of $v_m$ satisfies  $\varphi_{\psi}(v,v_m)
\leq 0$.

\end{defi}

\begin{lemma} The gradient lines of the elements of $\Phi$ are geodesics in
$\mathcal{K}$.

\end{lemma}

\proof First by \ref{lem;thirdpoint} we have  $|\varphi_{\psi} (w,v)| \leq \dist(w,v)$, for all $v$ and $w$. Now on a gradient line, we have by definition $\varphi_{\psi}(v_{n+1},v_n)=1$ for every consecutive vertices $v_{n}$ and $v_{n+1}$ on $l_\psi$. Therefore, as
$\varphi_{\psi}$ is a global cocycle, $\varphi_{\psi}(v_{m},v_n)=m-n$. The
triangular inequality gives the other inequality: $|\varphi_{\psi} (w,v)| \geq
\dist(w,v)$, and this proves the claim. $\square\smallskip$

We now state and prove the main property of the elements of
$\Phi$.

\begin{prop}(Coherence of gradient lines)\label{prop;coherence}

Let $\psi \in \Phi$. All its gradient lines are asymptotic to each
other. In other words they all converge to the same element of
$\partial \mathcal{K} \cup \mathcal{V}_{\infty}$.

\end{prop}

Before giving the proof of Proposition \ref{prop;coherence} we need some preliminary lemmas.

We say that two gradient lines of $\psi$ are \emph{divergent} if they are with different end points which are in the boundary or in the set of vertices of infinite valence.

\begin{lemma}\label{lemma;adjdivline}

Given two divergent gradient lines of $\psi$ one can find two divergent gradient lines
starting at the same vertex, or at two adjacent vertices.

\end{lemma} 

\proof Let $l_1$ and $l_2$ be two divergent gradient lines, and $v_1$ and
$v_2$ vertices on them. On a geodesic segment $[v_1, v_2]$,
consider $v$ the first vertex from which there is a gradient line
$l$ divergent from $l_1$. Either $v=v_1$ (and we are in the first
case of the lemma), or there is a vertex, $v'$, of $[v_1, v]$
adjacent to $v$. By definition of $v$, all gradient lines starting
at $v'$ are asymptotic to $l_1$, and we are in the second case of
the lemma. $\square\smallskip$

Let $l_1$, $l_2$ and $l_3$ be three geodesic lines connecting respectively $x_3,x_1$, $x_3,x_2$ and $x_1,x_2$ and suppose $l_1$ and $l_2$ are divergent gradient lines (possibly infinite) for $\psi$. 

\begin{lemma}\label{lem;center}

If $v$ is a vertex of $l_1$ such that $(x_1 \cdot x_2)_{x_3} - 50 \delta<\dist(x_3,v)< (x_1 \cdot x_2)_{x_3} + 50 \delta $ and if $\Ang_v(l_1)$ is more than $5\theta$ then the line $l_2$ pass through $v$ and $\Ang_v(l_1,l_2) \leq \theta$.

\end{lemma}

\proof It is enough to show that $v$ is on $l_2$, as we know by Lemma \ref{lem;thirdpoint} that two gradient lines starting at the same point make an angle at this point
smaller than $\theta$.

At distance $(x_1\cdot x_2)_{x_3} - 100 \delta$ from $x_3$ , we
connect $l_1$ and $l_2$ by a segment $\alpha_3$ of length less than
$10\delta$ , and we connect $l_1$ to $l_3$, and $l_3$ to $l_2$ at
distance $(x_1\cdot x_2)_{x_3} + 100 \delta$ from $x_3$ by two
others segments $\alpha_1$,$\alpha_2$ of length less than $10\delta$. In fact we will assume that such segments exist since that is the most complicated case, and otherwise one can reformulate the arguments in the next paragraph replacing the segments $\alpha_i$ by the vertices $x_i$ in order to obtain the same result, namely $v$ is either in $l_2$ or in $l_3$. 

Thus we have a loop of length less than $1000 \delta$ consisting of the 3 segments of lenght less than $10\delta$ and the subsegments of each $l_i$ lying in between the endpoints of these three segments. By Proposition \ref{prop;remarks}, no circuit of this length contains an angle more than $1000 \delta \leq \theta$. Hence, if $v$ is a vertex of $l_1$ such that $(x_1 \cdot
x_2)_{x_3} - 50 \delta<\dist(x_3,v)< (x_1 \cdot x_2)_{x_3} + 50
\delta $ and if $\Ang_v(l_1)$ is more than $5\theta$ then either
$l_2$ or $l_3$ pass through $v$ as the segments
connecting $l_1$, $l_2$ and $l_3$ are $10\delta$ short and
$50\delta$ far from $v$.

Now we assume that $v$ is on $l_3$ and not on $l_2$, and we now want to find a contradiction.
In this case, we have a
simple path from $x_3$ to $x_2$ consisting of the concatenation of
the piece of $l_1$ between $x_3$ and $v$ and the piece of $l_3$
between $v$ and $x_2$. The hyperbolicity of the space ensures that
this path remains at distance less than $60\delta$ from $l_2$. We have now two possiblities, either $l_2$ goes at infinity or $x_2$ is in $\mathcal{V}_{\infty}$.

In first case we
consider two adjacent vertices on $l_2$, $w$ and $w'$ such that
$\dist(v,w)<\dist(v,w')\leq 60\delta +1$ (two such vertices
necessarily exist since $l_2$ is a geodesic going to infinity). We choose a
 geodesic segment $[w',v]$ containing $w$.
If $\Ang_v([v,w'],[v,x_1]) < \theta$, then by triangular
inequality for angles, $\Ang_v([v,w'],[v,x_3]) \geq (5\theta-\theta)$,
and therefore Lemma \ref{lem;angles_triangles} for the
triangle $(x_3,v,w)$ implies that $v\in l_2$. On the other hand,
if $\Ang_v([v,w'],[v,x_1]) \geq \theta$,
let us subdivide   $[v,w']$ into maximal subsegments with maximal angle at most $\theta$.
Two consecutive such subsegments make an angle greater than $\theta$.
 Lemma \ref{lem;thirdpoint} (concerning the third point of  Proposition \ref{prop;locbusrad}) applied
successively to each of these subsegments, proves that
$[w',v]$ is a gradient line.
 This is a
contradiction, since the edge $(w,w')$ would be a gradient line in
both directions.

In the second case there is no vertex $w$ adjacent to $x_2$ with $\psi(w,x_2)=1$. Now again  if $\Ang_v([v,x_2],[v,x_1]) < \theta$, then by triangular
inequality for angles, $\Ang_v([v,x_2],[v,x_3]) \geq (5\theta-\theta)$,
and therefore Lemma \ref{lem;angles_triangles} for the
triangle $(x_3,v,x_2)$ implies that $v\in l_2$. On the other hand,
if $\Ang_v([v,x_2],[v,x_1]) \geq \theta$,
we subdivide  $[v,x_2]$ into maximal subsegments with maximal angle at most $\theta$ and again by
 Lemma \ref{lem;thirdpoint} applied
successively to each of these subsegments, we prove that
$[x_2,v]$ is a gradient line and hence the vertex $w$ adjecent to $x_2$ in $ [x_2,v]$
 satisfies $\psi(w,x_2)=1$, which is a contradiction. 

This completes to proof of Lemma \ref{lem;center}. $\square\smallskip$

\smallskip

\noindent\textit{Proof of Prop. \ref{prop;coherence} }. We argue by contradiction, and assume that $l_1$ and $l_2$ are two divergent gradient
lines. Now by Lemma \ref{lemma;adjdivline} we can assume that $l_1$ and $l_2$ are two divergent gradient
lines starting at the same vertex, or at two adjacent vertices.
In order to simplify the presentation we will only explicit the first case as the second case can be treated exactly with the same arguments.
Thus we assume that $l_1$ and $l_2$ start at the same vertex, hence by Proposition \ref{prop;visual}, there is a geodesic (possibly
bi-infinite) $l_3$, such that $(l_1,l_2,l_3)$ is a geodesic triangle with
vertices $x_1,x_2,x_3$ (see Figure 1), with $x_1$
and $x_2$ possibly at infinity.

>From the
previous lemma, and the definition of the Gromov product, we see
that a vertex in $l_1$ satisfying the assumption of the previous
lemma is in fact located at distance less than  $(x_1 \cdot
x_2)_{x_3}$ from $x_3$. Let $v$ be the last vertex satisfying the assumptions of the
previous lemma (or, if there is none, the vertex on $l_1$ at
distance $(x_1 \cdot x_2)_{x_3} - 50 \delta$ from $x_3$).
The two rays do not have an angle larger than $5\theta$ after $v$, until they arrive
at distance $50\delta$ from the small segments
connecting $l_1$, $l_3$, and $l_2$, $l_3$, because by the previous lemma, they would both pass by this vertex. Consider the cone
centred on the first edge of $l_1$ after our vertex $v$, of angle and radius
$10\theta$.

In this cone, let us parameterize the two subsegments of the two lines $l_i:[0,T_i]\to Cone_{10\theta, 10\theta} (v,e)$, for $i=1,2$.
 We know that $\dist(l_1(0),l_2(0))\leq  \delta$, because the triangle $(l_1,l_2,l_3)$ is $\delta$-thin.

  Moreover, as $\Ang_v(l_1,l_2) \leq \theta$ (by \ref{lem;center}) and as the two rays do not have angle larger than $5\theta$ after $v$, until they arrive close to the transition segments, the vertices $l_1(T_1)$ and $l_2(T_2)$ are at distance at least $(x_1 \cdot x_2)_{x_3} + 20 \delta$ from $x_3$, and therefore, at distance at least $10\delta$ from each other (see Figure 1), or possibly a segment $l_i$ reach $x_i$,
which, in this case, belong to the cone.

 By definition of our
alphabet $\mathcal{A}$, there must be a Busemann or a radial
cocycle whose restriction on this cone gives rise to the same
segments of gradient lines. This rules out the second case, and in
the first case, by hyperbolicity, two geodesic rays with such
subsegments would diverge at infinity, and we know that this
cannot happen for gradient lines of Busemann or radial cocycles.
This is a contradiction, and it proves the proposition. $\square\smallskip$

\begin{figure}[ht]\label{fig;triangle_and_cone}
  \begin{center}
    \includegraphics[height=3.5cm, width=7cm, angle=0]{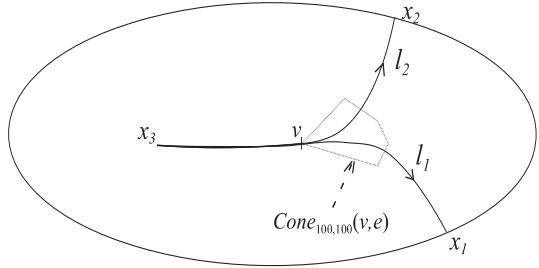}
  \end{center}
  \caption{Gradient lines and cone at the center of the triangle}
\end{figure}

We can now define the map $\Pi : \Phi \to \partial \Gamma$. For an
element $\psi$ in $\Phi$, we associate  $\Pi(\psi) \in
\partial \Gamma$, the point to which any gradient line of
$\varphi_{\psi}$ converge.

\subsection{End of the proof of Theorem \ref{theo;main}}

In order to complete the proof we need to show that $\Pi : \Phi
\to \partial \Gamma$ satisfies the Definition \ref{def;def2}
(Lemma \ref{lem;map}), and secondly that the action of $\Gamma$ on
$\partial \mathcal{K}$ is expansive (Proposition \ref{prop;exp}).

\begin{lemma}\label{lem;map}

The map $\Pi : \Phi \to \partial \Gamma$ is surjective continuous
and equivariant.

\end{lemma}

\proof Given a point $\xi$ in $\partial \Gamma$, by the visibility property, 
one can find a Busemann or a radial cocycle $\sigma$ associated to $\xi$. 
We consider the function $\tilde{\sigma} : G \to \mathcal{A'}$
obtained by $\tilde{\sigma}(\gamma)= \sigma|_{\gamma  Cone_{R,\Theta}(e_0,v_0)}$. Moreover, for every left coset $\gamma G_i$, $(i=1..m)$, we fix a representant $\tilde{\gamma}G_i$.
We consider the function, $\psi_i : G_i \to \mathcal{A}_i$, in $\Phi_i$  such that, if for some vertex $v$ and $g\in G_i$ with $\tilde{\gamma}g e_i=(\tilde{\gamma} p_i,v)$ one has $\tilde{\sigma} (\gamma^{-1}v,p_i)< 0$   then $\psi_i$ takes the value $\$ $ at some such $g$.    

For each $i$, the partition of $\Gamma$ in $G_i$-left cosets defines a
function $\tilde{\psi}_i : \Gamma \to \mathcal{A}_i$. 

This gives an function $\varphi:\Gamma \to \mathcal{A}$, namely the product of
$\tilde{\sigma}$ and of the $\tilde{\psi}_i$, $i=1,\dots,m$ which is, by Proposition 
\ref{prop;locbusrad}, an element of $\Phi$ whose gradient lines converge to
$\xi$. Thus, the map is surjective.

 If a sequence $\psi_n$ of elements of $\Phi$ converges to $\psi$, then, given a finite
set $E$ of edges of $\mathcal{K}$, the
gradient lines of $\psi_n$ will coincide with the gradient lines
of $\psi$ on $E$ for $n$ sufficiently large. If the point at infinity $\xi$ defined
by $\psi$ is in $\partial \mathcal{K}$, then by the usual topology of the
boundary of $\mathcal{K}$, the   
sequence of  points at infinity $\xi_n$, 
defined by the gradient lines of $\psi_n$, converges to $\xi$. If now $\xi$ is
in $\mathcal{V}_{\infty}$, then, by convergence criterion given in paragraph
2.4, the sequence $\xi_n$ converges to $\xi$. 
This ensures
the continuity of $\Pi$. 

Finally, if $\psi$ is an element of $\Phi$, and $\rho$ one of its gradient line, then $\gamma \rho$ is a gradient line for $\gamma \psi$. The point at infinity of $\gamma\rho$  is by definition $\gamma \Pi(\psi)$, and it is therefore equal to $\Pi(\gamma \psi)$.  Hence the map is equivariant. $\square\smallskip$

\begin{prop}(Expansivity)\label{prop;exp}

The action of a relatively hyperbolic group on its boundary is
expansive.

\end{prop}

\proof If $\Delta$ is the diagonal of $(\partial \Gamma)\times (\partial
\Gamma)$, then we have to find a neighborhood $U$ of $\Delta$ such
that $\Delta = \bigcap_{\gamma \in \Gamma} \gamma U$.

Let $\{e_1,\dots,e_m\}$ be a set of orbit representatives of the edges in
$\mathcal{K}$. Let $X$ be the set of pairs of points $(\xi_1,
\xi_2) \in (\partial \mathcal{K})^2$ such that there is a
bi-infinite geodesic between $\xi_1$ and $\xi_2$ passing through
one of the $e_i$. Let now  $\{ p_1,\dots,p_l \}$ be a set of orbit
representatives of the infinite valence points. Because they are
bounded parabolic points, the stabilizer $G_i$ of $p_i$ acts on
$\partial \Gamma \setminus \{ p_i \}$ with compact quotient. Let
then $Y$ be the set of pairs of points $(p_i,\zeta)$ where $\zeta$
is in a chosen compact fundamental domain for the action of $G_i$.

We now choose $U=(\partial \Gamma\times \partial \Gamma) \setminus
(X \cup Y)$. First we show that $\Delta = \bigcap_{\gamma \in
\Gamma} \gamma U$. The direct inclusion is trivial.

Let $(\xi_1,\xi_2) \in U$ and assume it is not in $\Delta$. We
will show that it is not in every translate of $U$. We consider two
cases, either $\xi_1, \xi_2$ are both in $\partial \mathcal{K}$, or one
them, say $\xi_1$, is a vertex of $\mathcal{K}$ of infinite
valence. In the first case, there is a bi-infinite geodesic from
one point to another, and it can be translated so that its image
passes by one of the $e_i$. Therefore, there is $\gamma$ such that
$\gamma(\xi_1,\xi_2)$ is in $X$, hence not in $U$. In the second
case, there is $\gamma \in \Gamma$ such that $\gamma \xi_1 $ is
one of the $p_i$. Now there is $\gamma' \in G_i$ such that
$\gamma'\gamma (\xi_1,\xi_2)$  is in $Y$, hence not in $U$. This
proves that the intersection of the translates of $U$ is equal to
the diagonal set.

Now we have to show that $U$ is a neighborhood of $\Delta$. That
is to say that a sequence of elements in $X\cup Y$ cannot converge
to a point of $\Delta$.

Let $(x_n = (\xi_1^n,\xi_2^n))_n$  be a converging sequence of
elements of $X$. After passing to a subsequence, one can assume
that, for all $n$, there is a bi-infinite geodesic between
$\xi_1^n$ and $\xi_2^n$ passing through a same edge $e_i$. If
$\xi_1^n \to \zeta_1$ and $\xi_2^n \to \zeta_2$, we see that
$\zeta_1$ and $\zeta_2$ are linked by a geodesic passing through
$e_i$, hence non-trivial. Therefore $\zeta_1 \neq \zeta_2$.

Let now $(y_n = (\xi_1^n,\xi_2^n))_n$  be a converging sequence of
elements of $Y$. After extraction, and without loss of generality,
one can assume that $\xi_1^n = p_i$, for all $n$, and for some
$i$. Then,  $\xi_2^n$ is in a compact fundamental domain for $G_i$
in $\partial \Gamma \setminus \{p_i\} $, and therefore does not
converge to $p_i$. This finally proves that $U$ is a neighborhood
of $\Delta$, and ends the proof of Proposition \ref{prop;exp}.
$\square\smallskip$

\section{Groups having the property of special symbol}\label{examples}

In this section we give examples of groups with the property of special symbol, and we
introduce a condition for it.

Let us begin with a necessary condition.

\begin{prop}\label{prop;finigene}

If $\Gamma$ has the property of special
symbol, then $\Gamma$ is finitely generated.

\end{prop}

\proof Let $\pi : \Phi \to \Gamma \cup \{\infty\}$ be a finite presentation with
special symbol. Let $\mathcal{A}$ be the alphabet. Let
$\mathcal{C}$ be a cylinder defining $\Phi$, and itself defined by
a  non-empty finite subset $F$ of $\Gamma$ and a set, $M$, of maps from $F$
to $\mathcal{A}$. The set of translates of $F$ is a covering of
$\Gamma$. Let $P$ be the nerve of the covering. As $F$ is finite,
$P$ is a finite dimensional, locally finite complex on which
$\Gamma$ acts properly discontinuously cocompactly. The set of
vertices of $P$ is naturally identified with $\Gamma$. The claim
is that $P$ is connected. If it was not, there would be distinct
connected components, $C_i$. Let $\gamma_i$ be a vertex of  $C_i$, and consider
$\sigma_i \in \Phi$ such that $\pi(\sigma_i) = \gamma_i$. Let
$\sigma \in \mathcal{A}^{\Gamma}$ such that $\sigma|_{C_i} \equiv
\sigma_i|_{C_i}$. Now, $\sigma$ has several special symbols (one
in each $C_i$). On the other hand all the cylinder conditions
defining $\Phi$ are satisfied, as by definition they are read on
the connected components of $P$. This is a contradiction, and it
proves the claim. Therefore, $\Gamma$ is generated by $F$ which is
a finite set. $\square\smallskip$

The next proposition is in fact is a slight variation
of a theorem of Gromov, a detailed proof of which can be
found in \cite{CP} (Corollary 8.2).

\begin{prop}\label{lem;hyperbolictrivialised}

If $\Gamma$ is a hyperbolic group, then its has the property of special symbol.

\end{prop}

\proof We repeat the proof of the main
theorem, seeing $\Gamma$ relatively hyperbolic relative to the
trivial subgroup $\{1\}$. A Cayley graph plays the role of
$\mathcal{K}$, and we consider the same cocycles. They can define
either a point at infinity, or a vertex of the graph. Thus, we
obtain our presentation choosing the special symbol to be the
restriction of a radial cocycle. $\square\smallskip$

Although it could be seen as a consequence of the proposition
above, the example 2 in part 1 already gave the basic examples of
$\mathbb{Z}$ and of finite groups. Most of our remaining examples
come from the following remarks.

\begin{prop}\label{prop;split}

If a group $\Gamma$ splits in a short exact sequence $\{1\} \to N
\to \Gamma \to H \to \{1\}$, and if both $N$ and $H$ have the property of special symbol,
then  $\Gamma$ has the property of special symbol.

\end{prop}

\begin{prop}\label{prop;virtual}

Let $G$ be a subgroup of finite index of a group $\Gamma$.
  The group $G$ has the property of special symbol if, and only if,  $\Gamma$ has the property of special symbol.

\end{prop}

Before giving the proofs, we give a consequence. A group $\Gamma$ is said
to be \emph{poly-hyperbolic} if  there is a sequence of subgroups $\{1\}=N_0
\lhd N_1 \lhd ... \lhd N_{k-1} \lhd N_k = \Gamma$, with all the quotients
$N_{i+1}/N_i$ hyperbolic.

\begin{coro}

 Every poly-hyperbolic group has the property of special symbol.
 In particular, this includes virtually polycyclic (hence, also virtually nilpotent) groups.

\end{coro}

\proof If $\Gamma$ is poly-hyperbolic, there is a sequence of subgroups
$\{1\}=N_0 \lhd N_1 \lhd ... \lhd N_{k-1} \lhd N_k = \Gamma$, with
all the quotients $N_{i+1}/N_i$ hyperbolic. Using the Proposition
\ref{prop;split}, and the fact that hyperbolic groups have the property of special symbol, an induction on $i$ tells that each  $N_i$
has the property of
special symbol, and especially $N_k$ which is $\Gamma$. $\square\smallskip$

\hskip .4cm

\noindent\textit{Proof of Prop. \ref{prop;split}.} Let us denote by $\mathcal{A}_N$, $\mathcal{A}_H$, $\$_N$, $\$_H$,
$\mathcal{C}_N$,  $\mathcal{C}_H$, $\Phi_N$, $\Phi_H$, the
alphabets, special symbols, cylinders, and subshifts of finite
type for the presentations of $N\cup \{\infty\}$ and $H \cup
\{\infty\}$. Let $F_N$, $F_H$, $M_N$ and $M_H$ be the finite
subsets of $N$ and $H$, and the sets of maps defining the two
given cylinders. From Proposition \ref{prop;finigene}, $N$ is
finitely generated, then up to enlarging $F_N$, we can assume that
$F_N$ generates $N$ (in fact, in the proof of Proposition
\ref{prop;finigene}, it is proved that necessarily, $F_N$
generates $N$). Let $\mathcal{A} =\mathcal{A}_H \times
\mathcal{A}_N$. Let us choose $\widetilde{H}$  a set of
representative of $H$ in $\Gamma$, and  for an element $h$ in $H$, 
we write $\tilde{h}$ for the element of $\widetilde{H}$ that 
maps on $h$ by the quotient map. Let $F$ be the finite
subset of $\Gamma$ defined by  $F=\{ \tilde{h}.n, \;  h \in F_H, n\in F_N \}$.
 Let $M$ be the following set
of maps: $M=\{ (m:F\to \mathcal{A}),  \exists m_H \in M_H, \forall
n\in F_N, m(\cdot.n)_1=m_H \; ; \;  \forall h\in F_H, 
m(\tilde{h}.\cdot)_2\in M_N \}$, where the subscripts $1$ and $2$
denote the coordinates in the product $\mathcal{A} = \mathcal{A}_H
\times \mathcal{A}_N$. Consider the cylinder defined by $F$ and
$M$, and the associated subshift of finite type, $\Phi$. We need
the following lemma.

\begin{lemma}

For any $\sigma \in \Phi$, there is at most one element $\gamma
\in \Gamma$ such that $\sigma(\gamma)=(\$_H,\$_N)$.

\end{lemma}

\proof We first prove that for any $\sigma \in \Phi$, there is at most
one left coset of $N$, $\tilde{h}N$, such that $\forall n \in N,
\; \sigma(\tilde{h}.n)_1 = \$_H$. By definition of $M$, if $n \in
F_N$, $n_0\in N$, then  $\sigma(\tilde{h}.n_0.n)_1$, the first
coordinate of $\sigma(\tilde{h}.n_0.n)$ only depends on
$\tilde{h}$ and $n_0$. But $F_N$ was chosen generating $N$, hence
$\sigma(\tilde{h}.n_0.n)_1$ only depends on $\sigma(\tilde{h})$.
But, by definition of $M$, the map $h\in H \mapsto
\sigma(\tilde{h})_1$ is in $\Phi_H$, and therefore, by the special
symbol property, there is at most one value of $\tilde{h}$ where
it takes the value $\$_H$, this proves the first step of the
lemma. 
Now, as the map  from $N$ to $ \mathcal{A}_N$ defined by 
$(n \mapsto \sigma(\tilde{h}.n)_2)$ is in $\Phi_N$,  if $\tilde{h}$ is such that
$\sigma(\tilde{h}.n)_1=\$_H $, there is at most one $n\in N$ such
that   $\sigma(\tilde{h}.n)_2=\$_N$. 
 This proves the lemma. $\square\smallskip$

Now, we define the map $\pi$ so that it sends a element $\sigma
\in \Phi$ on the point at infinity, if $\sigma$ does not contain
the symbol $(\$_H, \$_N)$, and on $\gamma \in \Gamma$ if
$\sigma(\gamma)=(\$_H, \$_N)$. The map $\pi$ is well defined, and
gives a finite presentation with special symbol of $\Gamma \cup
\{\infty\}$. $\square\smallskip$

\noindent \textit{Proof of Prop. \ref{prop;virtual}.} Assume  that $\Gamma$ has the property of
special symbol, and let $\mathcal{A}_\Gamma$, $\$_\Gamma$,
$\mathcal{C}_{\Gamma}$, $\Phi_{\Gamma}$ the alphabet, special symbol, cylinder,
and subshift of finite type associated. The cylinder is defined, as before, by
two sets : $F_\Gamma \subset \Gamma$ and $M_\Gamma \subset
\mathcal{A}^{F_\Gamma}$. We consider $\gamma_1,\dots,\gamma_n$ a set of orbit
representatives of left coset of $G$ in $\Gamma$, and we choose $
F=(\bigcup_{i=1}^n \gamma_i^{-1} F_\Gamma)\cap G$, a finite subset of $G$. We
set $\mathcal{A}= (\mathcal{A}_\Gamma)^n$ and $M\subset \mathcal{A}^F$ is the
set of the maps $m$ from $F$ to  $(\mathcal{A}_\Gamma)^n$ such that there exists
$m_\Gamma \in M_\Gamma$ whose translates $\gamma_i^{-1}m_\Gamma$ coincide with
the $i$-th coordinate of $m$. Those three choices define a  subshift of finite
type $\Phi \subset \mathcal{A}^G$. 
By definition of $M$, one sees that there is
a natural map $\Phi \to \Phi_\Gamma$ which consists of defining $\sigma_\Gamma (\gamma)$ by the ith coordinate of $\sigma (\gamma^{-1}_i \gamma$ if $\gamma$ is in the coset $\gamma_i G$. This map is a bijection, its inverse being the map
that associates to $\varphi\in\Phi_\Gamma$ the element $\sigma\in \Phi$ whose
$i$-th coordinate coincide with  $\gamma_i^{-1}\varphi$. Therefore, one has 
maps $ \Phi \to \Gamma \cup\{\infty\} \to G\cup\{\infty\}$, the second map being
identity on $G$ and sending each $\gamma_i$ to $1$. At this point we do not have
a special symbol, but, by property of $\Phi_\Gamma$, an element of $\Phi$ can
take a value in $\mathcal{A}$ which has $\$_\Gamma$ among its coordinates, only
once. Hence, by renaming each of those symbol by a single one $\$ $, we get the
expected presentation with special symbol.

Conversely, it suffices to see that the intersection of all the conjugates of
$G$ is  of finite index in $\Gamma$ (hence it has the property of special symbol). It is normal and of finite index in
$\Gamma$, and we can apply Proposition \ref{prop;split}. $\square\smallskip$


{\footnotesize



\thebibliography{99}
\bibitem{Bgrhyp}{\it B.H. Bowditch}, Notes on Gromov's hyperbolicity criterion for path-metric spaces, Group theory from geometrical viewpoint, ed. E.Ghys, a.Haefliger, A.Verjovsky, World Scientific, (1991), 64-167.
\bibitem{Brel}{\it B.H. Bowditch}, Relatively hyperbolic groups, preprint, University of Southampton, 1997 revised 1999.
\bibitem{Bbord}{\it B.H. Bowditch}, Boundaries of geometrically finite
  groups, Math Z. {\bf 230} (1999) n$^o$3, 509-527.
\bibitem{Bgf}{\it B.H. Bowditch}, Geometrically finiteness with variable negative curvature, Duke Math. J. {\bf 77} (1995) 229-274.
\bibitem{CP}{\it  M. Coornaert and A. Papadopoulos}, Symbolic dynamic and hyperbolic groups, Lecture Notes, Springer 1993.
\bibitem{CP2}{\it  M. Coornaert and A. Papadopoulos}, Horofunctions and
  symbolic dynamics on Gromov hyperbolic groups, Glasgow Math. J. {\bf 43}
  (2001), 425-456.
\bibitem{Dah}{\it  F. Dahmani}, Classifying spaces and boundaries for relatively
hyperbolic groups, Proc. London Math.Soc, {\bf 86} (2003), 666-684. 
\bibitem{Dcomb}{\it F.Dahmani}, Combination of convergence groups, preprint (2002). 
\bibitem{theseDah}{\it F.Dahmani}, Les groupes relativement hyperboliques et
  leurs bords, PhD thesis, IRMA Strasbourg, 2003. 
\bibitem{F}{\it B. Farb}, Relatively hyperbolic groups, Geom. Funct.
Anal. {\bf 8} (1998), n$^o$ 5, 810-840. \bibitem{Fri}{\it D. Fried}, Finitely
presented dynamical systems, Erg. th. \& Dyn. Syst. {\bf 7}, 1987, 489-507.
\bibitem{Gro}{\it M. Gromov},  Hyperbolic Groups, Essays in group theory, editor
S.Gersten, 75-263. 
\bibitem{Asli}{\it A. Yaman}, A topological characterisation
of relatively hyperbolic groups, preprint Southampton 2001 to appear in Journal f{\"u}r die reine und angewandte Mathematik (Crelle's Journal).
\bibitem{theseAs}{\it A. Yaman}, Boundaries of Relatively Hyperbolic Groups , PhD thesis, 2003.

}
\end{document}